\documentclass[english,a4paper,12pt]{article}
\usepackage{babel}
\usepackage[ansinew]{inputenc}
\usepackage{graphicx}
\usepackage{verbatim}
\usepackage{amsmath}
\usepackage{amssymb}
\def\C{\mathbb C}
\def\Cd{\widehat{\mathbb C}}

\newtheorem{theo}{Theorem}[section]
\newtheorem{ex}[theo]{Example}
\newtheorem{lemma}[theo]{Lemma}

\newtheorem{defi}[theo]{Definition}

\title{Meromorphic functions that partially share values with their first derivative}
\author{Andreas Sauer* \and Andreas Schweizer**}
\date{\empty}
\begin{document}
\maketitle
\thispagestyle{empty}
\textit{\small * Hochschule Ruhr West, Duisburger Str.~100, 45479 M\"{u}lheim an der Ruhr,\\ Germany, andreas.sauer@hs-ruhrwest.de}
\\[.3cm]
\textit{\small ** Department of Mathematics Education, Kongju National University, Gongju 32588,\\South Korea, schweizer@kongju.ac.kr}
\begin{abstract}
We consider uniqueness results for meromorphic functions $f \colon \C \to \Cd$ such that for certain values $a \in \C$ the implication $f(z)=a \Rightarrow f'(z)=a$ holds, i.e.~that $f$ and $f'$ share values {\it partially}. In particular, we give a result for four partially shared vales. 
\end{abstract}\vskip.3cm
\renewcommand{\thefootnote}{}
\footnotetext{\hspace*{-.51cm}2020 Mathematics Subject Classification: Primary: 30D35, Secondary: 30D45\\ %
Key words and phrases: Partially Shared Values, Uniqueness Problems, Picard value, Meromorphic Functions, Nevanlinna theory.}
\noindent
\section{Introduction}
Two meromorphic functions $f, g \colon \C \to \Cd$ are said to share the value $a\in\Cd$ if for all $z_0\in \C$ we have
$f(z_0)=a\Leftrightarrow g(z_0)=a$. There are famous and also beautiful theorems concerning shared values, the most famous being Nevanlinna's Five-Value Theorem, which says that if $f$ and $g$ share five values then $f \equiv g$ (e.g.~\cite{YaYi}, Chapter 3 and \cite{St}, section 3.3). Shared value problems are mainly part of Nevalinna's theory of meromorphic functions, we refer to \cite{ChYe}, \cite{St} and \cite{YaYi} (and the references given there).\\

Many authors have considered shared value problems for the special pair $f$ and $g=f'$ (see \cite{YaYi}, Chapter 8). Note that $f$ and $f'$ share $a=\infty$ automatically, so we will assume $a\in\C$ from now on. If $f$ and $f'$ share a value $a \in \C$ this can be considered as a restriction that resembles Hermite interpolation: Whenever at a point $z_0 \in \C$ one of the first two Taylor coefficients has the value $a$, then the other one is also $a$. The most famous theorem in this setting is the following result.
\\ \\
{\bf Theorem A.} \cite[Theorem 3]{Gu}, \cite[Satz 2]{MuSt}  \it
Let $f \colon \C \to \Cd$ be a non-constant meromorphic function that shares three finite 
values with $f'$, then $f\equiv f'$.
\rm
\\
\\
It is still an open problem, whether the number three in Theorem A can be replaced by two in case the shared values are non-zero.\\

In many proofs on shared values the given arguments do not always need the full equivalence $f=a \Leftrightarrow g=a$, but rather only $f=a \Rightarrow g=a$. It is natural to ask, which statements can be proved by using only such implications. This has led to the following definition.
\begin{defi} \rm A non-constant meromorphic function $f \colon \C \to \Cd$ shares the value $a\in\Cd$ {\it partially} with a non-constant meromorphic function $g \colon \C \to \Cd$, if for all $z_0\in \C$ we have
$$
f(z_0)=a\Rightarrow g(z_0)=a.
$$
\end{defi}
Note that the order in which the functions are named is important.\\

There is no uniqueness theorem for general $f$ and $g$ based only on partially shared values of $f$ and $g$. Simply choose any meromorphic function $f$, arbitrary many values $a_1, \ldots , a_n \in \C$ and let $E_k$ be the preimage of $a_k$ with respect to $f$. With the use of Weierstra{\ss} products and Mittag-Leffler series one can easily construct a meromorphic function $g$ such that $g(z)=a_k$ for every $z \in E_k$, so that all $a_k$ are partially shared by $f$ and $g$, but at the same time we are free to prescribe the values of $g'$ at all these points at will, e.g.~they could be chosen so large, that $f$ and $g$ are as different as we like.\\

As we will see, the situation for the pair $f$ and $f'$ is different. Partially shared values for entire functions $f$ and $f'$ have been treated in \cite{LuXuYi}, c.f.~Theorem B (see also \cite{SaSch}). Another instance where a partially shared value plays a role in connection with entire functions can be found in \cite{JaMuVo}, Satz 1. (See also \cite{YaYi}, section 8.1.5.) It is somewhat surprising, that we were not able to find uniqueness theorems concerning partially shared values of $f$ and $f'$ for general meromorphic functions, although the literature on partially shared values (and functions) in recent years, especially in connection with shifts $f(z+c)$ of $f$, is quite extensive (see e.g.~\cite{ChKoKu}). We could not even find a mentioning of the following theorem, which is easy to prove:

\begin{theo} \label{five} {\rm (Five Values Meromorphic)} \\
Let $f \colon \C \to \Cd$ be a non-constant meromorphic function that partially shares five finite values with its derivative, then $f \equiv f'$.
\end{theo}
\noindent
The same method of proof gives the following theorem for entire $f$.
\begin{theo} \label{three_entire} {\rm (Three Values Entire)} \\
Let $f \colon \C \to \Cd$ be a non-constant entire function that partially shares three finite values with its derivative, then $f \equiv f'$.
\end{theo}
Although the proofs of Theorem \ref{five} and Theorem \ref{three_entire} are straightforward, both theorems are sharp in the following sense: There exist meromorphic functions $f$ that partially share four finite values with $f'$ such that $f \not \equiv f'$ (see Theorem \ref{four} and Example \ref{example}, a)) and entire functions $f$ that partially share two finite values with $f'$ such that $f \not \equiv f'$ (see Theorem B).\\

\noindent
The following result in \cite{LuXuYi} improves on Theorem \ref{three_entire} considerably.\\ \\
{\bf Theorem B.} \cite[Theorem 1.2]{LuXuYi} {\rm (Two Values Entire)}\\
\it
Let $a$ and $b \neq 0$ be two distinct finite values, and let $f$ be a non-constant entire function. If $f$ and $f'$ partially share $a$ and $b$, then one of the following cases must occur:\\
(a) $f\equiv f'$,\\
(b) $f=A \exp\left( \frac{bz}{b-a} \right) + a$ with $a \neq 0$,\\
(c) $f=A \exp\left( \frac{az}{a-b} \right) + b$ with $a \neq 0$,\\
(d) $f=b \left( \frac{A}{2} \exp \left( \frac{z}{4} \right) + 1 \right)^2$ with $a=0$,\\
where $A$ is a non-zero constant.
\rm
\\ \\
\section{Results and Examples}
The main result of this paper is the following theorem.
\begin{theo} \label{four} {\rm (Four Values Meromorphic)} \\
Let $f$ be a non-constant meromorphic function that partially shares four finite values with its derivative.\\
a) If one of the four partially shared values is zero, then $f \equiv f'$.\\
b) If all four partially shared values are non-zero, then either $f \equiv f'$ or $f$ is of the form 
$$
f(z) = M(\text{\rm e}^{\lambda z})
$$
with $\lambda \neq 0$ and a M{\"o}bius-Transformation 
$$
M(w) = C \cdot \frac{w-A}{w-B}
$$ 
with non-zero constants $A$, $B$ and $C$, such that $M$ is non-constant ($A\neq B$) and such that
$$
\lambda \neq  \frac{\sqrt{B} - \sqrt{A}}{\sqrt{B} + \sqrt{A}} \quad \text{and} \quad \lambda \neq \frac{\sqrt{B} + \sqrt{A}}{\sqrt{B} - \sqrt{A}}.
$$
\end{theo}
\begin{ex} \label{example} \rm a) The suitable M{\"o}bius transformations $M$ in Theorem \ref{four} can be described in the following way: Given any M{\"o}bius transformation $M$ such that $a_1 := M(0)$ and $a_2 := M(\infty)$ are both non-zero and finite, then except for two values of $\lambda$ the meromorphic function $f(z) = M(\text{\rm e}^{\lambda z})$ partially shares four non-zero values $a_1, \ldots, a_4$ with its derivative. The values $a_1$ and $a_2$ are Picard exceptional values of $f$ and $a_3$ and $a_4$ are the values of $M$ at the two zeros of $\lambda w M'(w) - M(w)$ in $\C^*$. Note that the two values for $\lambda$ given in the theorem are independent of the (fixed) choice of $\sqrt{A}$ and $\sqrt{B}$. One such example is the function
$$
f(z) := \frac{4}{A \text{\rm e}^z + 1} - 3
$$
with $A \neq 0$. The Picard exceptional values of $f$ are $-3$ and $1$, and are therefore partially shared. The value $-1$ is even a shared value in the classical sense and $3$ is also a partially shared value.\\
\noindent
b) It is natural to ask whether the number four in Theorem \ref{four} can be reduced to three. The following three examples show that for three partially shared values the situation is more complicated than in Theorem \ref{four}. All examples we give here are of the form $R\left( \text{e}^{\lambda z} \right)$ with a rational function $R$ of degree two, which is in contrast to the conclusions in Theorem \ref{four}.

In the first example the three partially shared values are all non-zero. Set
$$
R(w) := \frac{18 w^2 + 18 w + 5}{18 w^2 + 36w + 20}.
$$
Then for $f(z) := R \left( \text{\rm e}^{z/6} \right)$, i.e.~with $w = \text{\rm e}^{z/6}$, the following can be shown by direct calculation: The values $1$, $1/4$ and $-1/2$ are partially shared values of $f$ and $f'$. The value $1$ is taken by $f$ if and only if $w = -5/6$, the value $1/4$ if and only if $w = -2/3$ and the value $-1/2$ if and only if $w = -2/3 \pm i/3$. For all these $w$ the function 
$$
w R'(w)/6 = \frac{3w(9w^2 + 15w + 5)}{2(9w^2 + 18w + 10)^2}
$$
takes the same values, which is equivalent to the statement that the three values are partially shared by $f$ and $f'$.

An example with three partially shared values, where one shared value is $0$, is given by
$$
R(w) := \left( \frac{w-2}{w-1} \right)^2
$$
with $w = \text{\rm e}^{-z/12}$. Obviously $0$ is a shared value of $R$ and
$$
 - w R'(w)/12 =  \frac{- w(w-2)}{6(w-1)^3}
$$
in $\C^*$. Further the value $1$ is taken by $f$ if and only if $w=3/2$ and the value $4$ if and only if $w=4/3$. At these two points $- w R'(w)/12$ also takes the values $1$ and $4$, respectively, so that $f(z) = R \left(\text{\rm e}^{-z/12}\right)$ and $f'$ share the value $0$ and partially share $1$ and $4$.

In the last example one of the partially shared values is a non-zero Picard exceptional value. Set
$$
R(w) = \frac{(A w + 2)^2}{A^2 w^2 + 4Aw + 8}
$$
with $A\neq 0$ and $f(z) := R \left( \text{\rm e}^{z/2} \right)$, i.e.~$w = \text{\rm e}^{z/2}$. Then $f$ has the Picard exceptional value $1$ which is therefore partially shared by $f$ and $f'$. Obviously $0$ is a shared value of $R$ and
$$
 w R'(w)/2 =  \frac{4Aw(Aw + 2)}{(A^2 w^2 + 4Aw + 8)^2}
$$
in $\C^*$. Further the value $1/2$ is taken by $f$ if and only if $w=-4/A$. It is easy to check that $w R'(w)/2$ also equals $1/2$ at this point. Hence $1$, $0$ and $1/2$ are partially shared by $f$ and $f'$, where $0$ is even a shared value in the classical sense.
\end{ex}
\noindent
If $f$ has only multiple poles we can improve the number four in Theorem \ref{four} to three.
\begin{theo} \label{multiple_poles} {\rm (Three Values Meromorphic with Multiple Poles)} \\
Let $f \colon \C \to \Cd$ be a non-constant meromorphic function that has only multiple poles and that partially shares three finite values with its derivative.\\
a) If the three partially shared values are non-zero, then $f \equiv f'$.\\
b) If one of the three partially shared values is zero, then either $f \equiv f'$ or $f$ is of the form 
$$
f(z) = \left( M \left(\text{\rm e}^{\lambda z} \right) \right)^2
$$
with a M{\"o}bius-Transformation 
$$
M(w) = C \cdot \frac{w-A}{w-B}
$$ 
with non-zero constants $A$, $B$ and $C$, such that $M$ is non-constant, $A^2 \neq B^2$ and
$$
\lambda = \frac{B-A}{4(A+B)},
$$
when the partially shared values are $0$, $C^2$ and $C^2 A^2/B^2$.
\end{theo}
\section{Proofs for the theorems}\label{proofs}

The following is a consequence of the lemma on the logarithmic derivative and basic properties of the logarithmic derivative concerning its poles. See e.g.~\cite{MuSt}, inequality (2) and e.g.~Theorem 3.5.1 in \cite{ChYe} for the statement on the error term if $f$ has finite order.

\begin{lemma}  \label{N(f'-f)} Let $f \colon \C \to \Cd$ be a non-constant meromorphic function with $f \not \equiv f'$, then
$$
N(r,0,f'- f) \le T(r,f) + \overline{N}(r,f) + S(r,f)
$$
and
$$
T(r,f'/f) \le \overline{N}(r,0,f) + \overline{N}(r,f) + S(r,f).
$$
If $f$ has finite order then $S(r,f) = O(\log(r))$.
\end{lemma}
\noindent
{\bf Proof of Theorem \ref{five} and Theorem \ref{three_entire}.} First let $f$ be meromorphic and $a_1, \ldots , a_5$ be the partially shared values. If all these values are non-zero, then the second fundamental theorem of Nevanlinna theory and Lemma \ref{N(f'-f)} show
\begin{align*}
3T(r,f) & \le \sum_{k=1}^5 \overline{N}(r,a_k,f) + S(r,f)\\ 
& \le N(r,0,f'- f) + S(r,f)\\
& \le T(r,f) + \overline{N}(r,f) + S(r,f)\\
& \le 2 T(r,f) + S(r,f)
\end{align*}
Dividing this inequality by $T(r,f)$ and letting $r \to \infty$ (outside a possible exceptional set for $r$) we get the contradiction $3 \le 2$.\\

If one of the values is $0$, say $a_5=0$, then $f$ only has multiple zeros, hence $\overline{N}(r,0,f) \le \frac{1}{2} N(r,0,f)$. Again the second fundamental theorem of Nevanlinna theory and Lemma \ref{N(f'-f)} show
\begin{align*}
2T(r,f) & \le \sum_{k=1}^4 \overline{N}(r,a_k,f) + S(r,f)\\ 
& \le N(r,1,f'/f) + S(r,f)\\
& \le T(r,f'/f) + S(r,f)\\
& \le \overline{N}(r,0,f) + \overline{N}(r,f) + S(r,f)\\
& \le \frac{1}{2} N(r,0,f) + \overline{N}(r,f) + S(r,f) \\
& \le \frac{3}{2} T(r,f) + S(r,f)
\end{align*}
Just as above this gives the contradiction $2 \le 3/2$.\\

For entire functions these arguments have to be adjusted only slightly. We omit the details. $\square$\\

\noindent
A key feature in the following proofs is that, under our assumptions, $f$ turns out to be a normal function.
 
\begin{defi} \rm A non-constant meromorphic function $f \colon \C \to \Cd$ is called {\it normal} if the spherical derivative $f^{\#}(z) = |f'(z)|/(1 + |f(z)|^2)$ is bounded in the whole plane.
\end{defi}
\noindent
As usual, we denote by $\rho(f)$ the {\it order} of growth of $f$ (\cite{St}, section 2.1.7). It is a well-known fact that $\rho(f) \le 2$ for normal functions. This can be proved by simple estimation of $T(r,f)$ in the form named after Ahlfors-Shimizu (\cite{ChYe}, Theorem 1.11.3 or \cite{St}, section 2.4.2).\\

The normality of $f$ under our assumptions follows directly from Theorem 2 in \cite{XuFa} (see also Corollary 1.3 in \cite{Sch}):

\begin{theo} {\rm (\cite{XuFa}, Theorem 2)} \label{normal} If a non-constant meromorphic function $f \colon \C \to \Cd$ partially shares three finite values with its derivative, then $f$ is normal. 
\end{theo}

\noindent
The following theorem (\cite{FrDu}, Th\'{e}or\`{e}me 1) is of crucial importance for our main proof.

\begin{theo} {\rm (\cite{FrDu}, Th\'{e}or\`{e}me 1)} \label{not_deficient} Let $f \colon \C \to \Cd$ be a non-constant normal meromorphic function and $a \in \Cd$, then
$$
T(r,f) = N(r,a,f) + O(r).
$$ 
\end{theo}

\noindent
Further we need the following lemma which is a direct consequence of a well-known version of the lemma on the logarithmic derivative, given e.g.~in Theorem 3.5.1 in \cite{ChYe}.

\begin{lemma} {\rm (Ngoan-Ostrovskii)} \label{Ngoan-Ostrovskii} Let $f \colon \C \to \Cd$ be a non-constant meromorphic function with order $\rho(f) \le 1$, then
$$
m(r,f'/f) = o(\log(r)).
$$ 
\end{lemma}

\noindent
{\bf Proof of Theorem \ref{four}.} Let $a_1, \ldots ,a_4$ be the partially shared values.\\

First we show that $f$ cannot be a rational function. Suppose $f$ is rational of degree $d>0$. We consider $Q := f'/f$. Since $Q$ has poles of first order only in the zeros and poles of $f$ it follows $\deg(Q) \le 2d$. If the values $a_1, \ldots ,a_4$ are all non-zero, then for at least three of them all preimages lie in $\C$ and are simple, so there are at least $3d$ of them. Since these points are $1$-points of $Q$ we get a contradiction. Now suppose $a_4=0$. Since all zeros of $f$ are at least of order two, we get by counting poles $\deg(Q) \le 3d/2$. For at least two of the values $a_1, \ldots ,a_3$ all preimages lie in $\C$ and are simple, so there are at least $2d$ of them. Just as above we get the contradiction $2d \le 3d/2$.\\

From now on we assume that $f$ is transcendental and $f \not \equiv f'$. Theorem \ref{normal} shows that $f$ is normal. We prove $\rho(f) \le 1$. Three of the partially shared values are non-zero, say $a_1, \ldots, a_3$. Then we get with Lemma \ref{N(f'-f)}:
\begin{equation} \label{sum}
\sum_{k=1}^3 N(r,a_k,f) = \sum_{k=1}^3 \overline{N}(r,a_k,f) \le N(r,f'- f) \le T(r,f) + \overline{N}(r,f) + O(\log(r)).
\end{equation}
Suppose $\rho(f) > 1$. Then there is a sequence $r_n \to \infty$, such that $r_n/T(r_n,f) \to 0$ (and of course $\log(r_n)/T(r_n,f) \to 0$). Divide (\ref{sum}) by $T(r,f)$ and set $r = r_n$. Then Theorem \ref{not_deficient} and taking limits give the contradiction $3 \le 2$.\\
We consider
$$
\Phi := \frac{f'(f'-f)}{(f-a_1)(f-a_2)(f-a_3)(f-a_4)}.
$$
$\Phi$ has no poles in poles of $f$, and $\Phi$ has no poles if $f=a_k$, even if $a_k=0$. Hence $\Phi$ is entire. By partial fraction decomposition $\Phi$ can be written as a combination of logarithmic derivatives of $f-a_k$. Since $\rho(f) \le 1$ Lemma \ref{Ngoan-Ostrovskii} shows $T(r,\Phi) = o(\log(r))$. This implies that $\Phi$ is constant. If $\Phi \equiv 0$ then the claim follows immediately. Therefore we suppose
$$
f'(f'-f) - c \cdot (f-a_1)(f-a_2)(f-a_3)(f-a_4) = 0
$$
with $c \neq 0$.\\

This gives a Briot-Bouquet equation $P(f,f')=0$ with
$$
P(x,y) = y(y-x) - c \cdot (x-a_1)(x-a_2)(x-a_3)(x-a_4).
$$
According to Theorem 5.12 in \cite{St} solutions of this equation are either rational functions, elliptic functions or functions of the form $R(\exp(\lambda z))$ with a rational function $R$. Since elliptic functions have order $2$ and $f$ cannot be rational we conclude $f(z) = R(\exp(\lambda z)$. Since the degree of $P$ with respect to $y$ is two, Theorem 5.12 in \cite{St} shows that $d := \deg(R) \le 2$ when $\lambda$ is chosen such that $\exp(\lambda z)$ has the same fundamental period as $f$, which we will assume from now on. The shared value problem for $f$ and $f'$ can be formulated as a shared value problem for $R(w)$ and 

$$
\partial_z R(w) := \lambda w R'(w).
$$
$f$ and $f'$ share the value $a$ partially in $\C$ if and only if $R$ and $\partial_z R$ share the value $a$ partially in $\C^*$. The assumption $f \not \equiv f'$ translates to $R \not \equiv \partial_z R$.\\

It is a basic property of $\partial_z R$ that if $R$ has a pole of order $k$ in $w=0$ or $w=\infty$, then $\partial_z R$ also has a pole of order $k$ in $w=0$ or $w=\infty$. If $R$ has a finite value of order $k$ in $w=0$ or $w=\infty$, then $\partial_z R$ has a zero of order $k$ in $w=0$ or $w=\infty$. This can easily be shown by Laurent expansion in $w=0$ and $w=\infty$.\\
We translate $\Phi$ into this notation:
$$
\Phi := \frac{\partial_z R \cdot (\partial_z R-R)}{(R-a_1)(R-a_2)(R-a_3)(R-a_4)}.
$$
Still we have $\Phi \equiv c \neq 0$, but now $\Phi(w)$, $R(w)$ and $\partial_z R(w)$ are defined for all $w \in \Cd$, in particular in $w=0$ and $w=\infty$. Let $n_{\infty}$ be the number of poles of $R$ in $\C^*$ counted without multiplicities. Counting poles it is easy to see that $\deg(\partial_z R) = d + n_{\infty}$ and $\deg(\partial_z R - R) \le d + n_{\infty}$.\\

To prove Theorem \ref{four}, a) we assume $a_4=0$. First we consider the case that $R$ has a zero in $\C^*$. Since this zero has to be double we have $d=2$ and hence $\deg(\partial_z R - R) \le 4$. If $R$ has simple points at $w=0$ and $w=\infty$, then $R$ has at least four preimages for $a_1, \ldots , a_3$ in $\C^*$. Together with the common zero of $R$ and $\partial_z R$ this gives at least five zeros for $\partial_z R - R$, a contradiction. Hence we can assume without loss of generality that the only preimage of $a_3$ lies in $w=0$ or $w=\infty$. Consider
$$
\Theta := \frac{\partial_z R}{R (R-a_3)}.
$$
$\Theta$ has no poles in poles of $R$ and no pole in the only $a_3$-point. The only simple pole of $\Theta$ lies in the common zero of $R$ and $\partial_z R$, thus $\deg(\Theta)=1$. For at least one of the values $a_1$ or $a_2$ there are two simple preimages in $\C^*$, say for $a_2$. In both of these points we have $\Theta = 1/(a_2 - a_3)$, a contradiction.\\
We conclude that $R$ has no zero in $\C^*$ so that $f$ has $0$ as a Picard exceptional value. Consider $h(z) := 1/f(-z)$. Then $h$ is an entire function with 
$$
h=1/a_k \Rightarrow h'= 1/a_k
$$
for $k=1,2,3$. Theorem \ref{three_entire} shows $h \equiv h'$ which implies $f \equiv f'$.\\
   
To prove Theorem \ref{four}, b) we assume that the values $a_1, \ldots ,a_4$ are all non-zero. It follows that for two of these values $w=0$ and $w = \infty$ have to be the only preimages (so that $f$ has these values as Picard exceptional values), since otherwise $\partial_z R - R$ has more than $2d$ zeros which contradicts $\deg(\partial_z R - R) \le d + n_{\infty} \le 2d$. Say $R(0) = a_1$ and $R(\infty)=a_2$ and $R$ has no $a_1$- or $a_2$-points in $\C^*$. Consider
$$
\Psi := \frac{\partial_z R}{(R-a_1)(R-a_2)}.
$$
From the above mentioned properties of $\partial_z R$ and the other assumptions it follows that $\Psi$ has no poles in $\Cd$, hence $\Psi \equiv C$ with a constant $C$. We get $C \neq 0$ since $\Psi$ is a factor of $\Phi$. This gives the seperable Ricatti equation 
$$
f' = C(f-a_1)(f-a_2)
$$
which can be solved elementarily. All solutions turn out to be of the form $f=M(\exp(C(a_2-a_1)z))$ with a M{\"o}bius transformation $M$. The rest of the proof of Theorem \ref{four}, b) is completely elementary. One only has to choose $M$ in such a way that $M(0)$ and $M(\infty)$ are both in $\C^*$, and choose $\lambda$ such that $\partial_z M - M$ has two simple zeros in $\C^*$. Direct calculation shows that $\partial_z M - M$ has only one zero of order two exactly for the two exceptional values of $\lambda$ given in the theorem. $\square$
\vskip.3cm
\noindent
{\bf Proof of Theorem \ref{multiple_poles}.} The main part of the proof runs along the lines of the proof of Theorem \ref{four} with only minor adjustments:\\

We exclude rational functions. Suppose $f$ is rational with degree $d$. First we assume that the partially shared values are all non-zero. Since $f$ has only multiple poles we get by counting poles that $Q := f'/f$ has degree at most $3d/2$ but at least $2d$ $1$-points in the preimages of $a_1, \ldots, a_3$, so that $Q$ is constant, a contradiction. Now we suppose $a_3=0$, so that all zeros of $f$ are multiple. Counting poles again we see that $Q$ has degree at most $d$. In the preimages of $a_1$ and $a_2$ the function $Q$ has $1$-points, so that one of the two values has to be omitted by $f$ in $\C$, say $f(\infty)=a_1$ and $f$ has no $a_1$-point in $\C$. Consider $g(z) = 1/(f(z)-a_1)$. Then $g$ is a polynomial of degree $d$. In the $d$ preimages of $a_2$ the derivative $g'(z) = -f'(z)/(f(z)-a_1)^2$ takes the value $-a_2/(a_2 - a_1)^2$. Since $\deg(g')=d-1$ it follows that $g' \equiv -a_2/(a_2 - a_1)^2$ and therefore $f$ is a M{\"o}bius transformation. Since $f(\infty)=a_1$ we conclude that $f$ has a simple pole in $\C$, contradicting the assumptions.\\

From now on we assume that $f$ is transcendental and $f \not \equiv f'$. Theorem \ref{normal} shows that $f$ is normal. We prove $\rho(f) \le 1$. Two of the partially shared values are non-zero, say $a_1, a_2$. Since all poles of $f$ are multiple we have $\overline{N}(r,f) \le \frac{1}{2} N(r,f)$.Then we get with Lemma \ref{N(f'-f)}:
\begin{align*} \label{sum_2}
 N(r,a_1,f) +  N(r,a_2,f) & = \overline{N}(r,a_1,f) +  \overline{N}(r,a_2,f)\\ 
& \le N(r,f'- f) \\
& \le T(r,f) + \overline{N}(r,f) + O(\log(r))\\
& \le T(r,f) + \frac{1}{2} N(r,f) + O(\log(r))\\
& \le \frac{3}{2} T(r,f) + O(\log(r)).
\end{align*}
For $\rho(f) > 1$ this leads by the same argument as in the forgoing proof to the contradiction $2 \le 3/2$.\\

We consider
$$
\Phi := \frac{f'(f'-f)}{(f-a_1)(f-a_2)(f-a_3)}.
$$
Since all poles of $f$ are multiple, $\Phi$ has no poles in poles of $f$. Again $\Phi$ is entire and can be written as a combination of logarithmic derivatives so that $\Phi$ is constant because of $\rho(f) \le 1$. Again $\Phi \equiv 0$ or
$$
f'(f'-f) - c \cdot (f-a_1)(f-a_2)(f-a_3) = 0
$$
with $c \neq 0$. Since we have excluded rational functions $f$ and $\rho(f) > 1$ we get with Theorem 5.12 in \cite{St} that $f(z) = R(\exp(\lambda z))$ with a rational function $R$ with $\deg(R) \le 2$.\\

\noindent
Again we consider $R$ and $\partial_z R$ on $\C^*$. We treat several cases:\\

1) If $R$ has no pole in $\C^*$ then $f$ is entire and Theorem \ref{three_entire} shows $f \equiv f'$.\\

2) If $R$ has a pole in $\C^*$ then it has to be exactly one double pole and $d=2$.\\

2a) If the three partially shared values $a_1,\ldots,a_3$ are all non-zero, then $h(z) = 1/f(-z)$ partially shares $1/a_1,\ldots,1/a_3$ and $0$ with $h'$. Theorem \ref{four}, a) shows $h \equiv h'$ and hence $f \equiv f'$.\\

2b) Suppose $a_3=0$.\\

2bi) If $R$ has no zero in $\C^*$, then $h(z) = 1/f(-z)$ is entire and partially shares the values $1/a_1$, $1/a_2$ and $0$ with $h'$. Theorem \ref{three_entire} shows $h \equiv h'$ and hence $f \equiv f'$.\\

2bii) If $R$ has a zero in $\C^*$, it has to be a double zero since $0$ is partially shared. It follows $R(w)=(M(w))^2$ with a M{\"o}bius transformation 
$$
M(w) = C \cdot \frac{w-A}{w-B}
$$ 
with non-zero constants $A$, $B$ and $C$, such that $M$ is non-constant ($A\neq B$). First note that the total ramification of $2$ for $R(w)=(M(w))^2$ lies in the double zero and the double pole. Hence all other values of $R$ lie in $\C^*$ and are simple. We have $\deg(\partial_z R)=3$ and $\deg(\partial_z R - R)=3$. If one of the values $R(0)$ or $R(\infty)$ is not partially shared, then we have at least three preimages for $a_1$ and $a_2$ in $\C^*$. Together with the common zero this gives four zeros for $\partial_z R - R$, a contradiction. First we exclude the possibility $R(0)=R(\infty)$. Suppose to the contrary $a_1 := R(0)=R(\infty)$ which is equivalent to $A^2=B^2$. Consider
$$
\Theta := \frac{\partial_z R}{R (R-a_1)}.
$$
Then $\Theta$ has only one simple pole in $w=A$, hence $\deg(\Theta)=1$. But at the two preimages of $a_2$ in $\C^*$ we have $\Theta = 1/(a_2 - a_1)$, a contradiction. Hence $a_1 := R(0) = C^2A^2/B^2$ and $a_2 := R(\infty)=C^2$ are distinct. The only finite solution of the equation $R(w) = C^2$ is $w=(A+B)/2 \in \C^*$. The requirement $\partial_z R((A+B)/2) = C^2$ leads with simple calculations to the unique choice
$$
\lambda = \frac{B-A}{4(A+B)}.
$$
Given that $\lambda$ it is easy to check that $a_1$ is also partially shared at $w=2AB/(A+B)$. $\square$

\noindent
{\bf Acknowledgements.} We would like to thank Alexandre Eremenko for bringing \cite{FrDu} to our attention, which contains an important argument in the proofs of our main results.

\end{document}